\documentclass[a4paper,11pt]{article}
\usepackage{graphicx}
\usepackage{multirow}
\usepackage{color}
\usepackage{latexsym}
\usepackage{amssymb}
\usepackage{amsfonts}
\usepackage{mathrsfs}
\usepackage{amsmath}
\usepackage{indentfirst}
\newcommand{\cvd}{\hfill $\blacksquare$\bigskip}
\newtheorem{definition}{Definition}[section]
\newtheorem{theorem}{Theorem}[section]
\newtheorem{proposition}{Proposition}[section]

\author{Elena Barcucci\thanks{Dipartimento di Matematica e Informatica ``U. Dini'', Universit\`a degli
Studi di Firenze, Viale
 G.B. Morgagni 65, 50134 Firenze, Italy. {
 \tt \ elena.barcucci@unifi.it,\quad antonio.bernini@unifi.it,\quad
renzo.pinzani@unifi.it}}\and Antonio Bernini$^*$ \and Renzo Pinzani$^*$}

\title{Strings from linear recurrences and permutations:
a Gray code}
\date{}
\begin{document}

\maketitle

\begin{abstract}
Each positive increasing integer sequence $\{a_n\}_{n\geq 0}$ can serve as a numeration system to represent each non-negative integer by means of suitable coefficient strings. We analyse the case of $k$-generalized Fibonacci sequences leading to the binary strings 
avoiding $1^k$. We prove a bijection between the set 
of strings
of length $n$ and the set of permutations of $S_{n+1}(321,312,23\ldots(k+1)1)$. Finally, basing on a known Gray code for
those strings, we define a Gray code for $S_{n+1}(321,312,23\ldots(k+1)1)$, where two consecutive permutations differ by an adjacent transposition. 
\end{abstract}


\section{Introduction}\label{intro}
In \cite{BSS} the authors asked for a combinatorial interpretation of 
the recurrence $f_{m+1}=6f_m-f_{m-1}$, with $f_0 = 1, f_1 = 7$ (sequence M4423 of \cite{S}). A general solution appeared in \cite{BR} where a combinatorial
interpretation for the recurrences of the form $a_m = ka_{m-1} + ha_{m-2}$ and underlying some conditions on $h$ and
$k$, was given. In particular, by considering the sequence arising from the recurrence as a numeration system \cite{F} able to represent each non-negative integer as strings (see next section), it is possible to completely characterize the language of these strings and give a recursive construction for such a language.

Recently (\cite{BBP1,BBP2}), a Gray code has been defined 
for the languages deriving from the recurrences $a_m = ka_{m-1} + 
ha_{m-2}$ for $k\geq h\geq 0$ and $a_m = ka_{m-1} - ha_{m-2}$ with
$k>h>0$ and $h$ even.
In this paper we continue this research considering the well known  $k$-generalized Fibonacci sequence. A first result shows that in this case the language is the set $F^{(k)}$ of the binary strings avoiding $k$ consecutive $1$'s (i.e. the pattern $1^k$) and vice versa: each binary string avoiding $1^k$ is the representation of a unique non-negative integer. Moreover, by reading  the strings of the set $F_n^{(k)}$ of the binary strings of $F^{(k)}$ of length $n$ as the inversion arrays of permutations of length $n+1$, we show that the strings are in bijection with the set $S_{n+1}(321,312,23\ldots(k+1)1)$ of avoiding $321$, $312$, and $23\ldots (k+1)1$ permutations \cite{B}. It is already known that the strings in $F_n^{(k)}$ can be listed in a Gray code order with Hamming distance \cite{H} equal to $1$ (see \cite{V}). We show that this Gray code \cite{G} can be transferred to $S_{n+1}(321,312,23\ldots(k+1)1)$ where two consecutive permutations differ only for an adjacent transposition (i.e. switching two consecutive entries). 

\section{Preliminaries}

Given a sequence $\{a_m\}_{m\geq 0}$ of integers such
that $a_0=1$ and $a_m<a_{m+1}$ for each $m\in \mathbb N$, let $N$ be any non-negative integer. Consider the largest term $a_n$ of the sequence such that $a_n\leq N$. More precisely, $a_n=\max\{a_m\ |\ a_m\leq N\}$ (for the particular case $N=0$, see below). We divide $N$ by $a_n$ obtaining $N=d_na_n+r_n$. Obviously, for the remainder $r_n$, it is clear that $r_n<a_n$. If we divide $r_n$ by $a_{n-1}$, we get $r_n=d_{n-1}a_{n-1}+r_{n-1}$, with $r_{n-1}<a_{n-1}$. Then, iterating this procedure until the division by $a_0=1$ (where of course the remainder is 0), we have:
$$
\begin{array}{rllll}
	N&=&d_na_n+r_n\quad &\phantom{iopppppp} & 0\leq r_n<a_n\ ,\\
	&&&&\\
	r_n&=&d_{n-1}a_{n-1}+r_{n-1} & & 0\leq r_{n-1}<a_{n-1}\ ,\\
	&&&&\\
	r_{n-1}&=&d_{n-2}a_{n-2}+r_{n-2} & & 0\leq r_{n-2}<a_{n-2}\ ,\\
	&&&&\\
	\cdots&=&\cdots\cdots\cdots\cdots&&\cdots\cdots\cdots\cdots\\
	&&&&\\
	\cdots&=&\cdots\cdots\cdots\cdots&&\cdots\cdots\cdots\cdots\\
	&&&&\\
	r_3&=& d_2a_2+r_2&&0\leq r_2<a_2\ ,\\
	&&&&\\
	r_2&=& d_1a_1+r_1&&0\leq r_1<a_1\ ,\\
	&&&&\\
	r_1&=& d_0a_0\ .\\
\end{array}
$$

The above relations imply that:
\begin{equation}\label{rappr}
	N=d_na_n+d_{n-1}a_{n-1}+d_{n-2}a_{n-2}+\ldots\ldots+d_1a_1+d_0a_0\ .
\end{equation}

\medskip
Expression (\ref{rappr}) is the representation of $N$ in the numeration system $S=\{a_0,a_1,a_2,\ldots\ldots\}$, and the string $d_nd_{n-1}\ldots d_1d_0$ is associated to the number $N$
(in what follows the term ``representation" equivalently refers either to 
expression (\ref{rappr}) or to its associated string).
This method \cite{F} can be applied to every non-negative integer and in the case 
$N=0$, clearly, all the coefficients $d_i$ are $0$
(in other words the representation of $0$ is simply the string $0$). Moreover, we have

\begin{equation}\label{resto}
	r_{i}=d_{i-1}a_{i-1}+d_{i-2}a_{i-2}+\ldots\ldots+d_1a_1+d_0a_0<a_i\ ,
\end{equation}

\noindent
for each $i\geq 0$.

\medskip
\noindent
It is possible to show \cite{F} that if $N=\sum_{i\geq 0}^{n}d_ia_i$ with

\begin{eqnarray}\label{minore}
	d_ia_i+d_{i-1}a_{i-1}+\ldots+d_1a_1+d_0a_0<a_{i+1}
\end{eqnarray}
\noindent
for each $i\geq 0$, then the representation $N=\sum_{i\geq 0}^{n}d_ia_i$ is unique. For the sake of completeness, we recall the complete theorem:

\begin{theorem}\label{teo}
	Let $1 = a_0 < a_1 < a_2 < \ldots$ be any finite or infinite sequence of integers. Any non-negative integer $N$ has precisely one representation in the system $S = \{ a_0, a_1, a_2, ... \}$ of the form $N=\displaystyle \sum_{i\geq 0}^{n}d_ia_i$ where the $d_i$ are non-negative integers satisfying (\ref{minore}).	
\end{theorem}

As an example, consider the well-known sequence of Pell numbers (sequence M1413 in \cite{S}) $p_m=1,2,5,12,29,\ldots$ defined by $p_0=1$, $p_1=2$, $p_m=2p_{m-1}+p_{m-2}$. The representation of $N=16$ is associated to the string $1020$.

\section{Strings from a number sequence}

Given a sequence $\{a_m\}_{m\geq0}$, for a fixed $m>0$, we consider all the integers $\ell\in \{ 0, 1, 
2,\ldots, a_m-1\}$. According to the scheme of the previous section, the 
representations of the integers $j$ with $a_{m-1}\leq j<a_m$ is 
$j=d_{m-1}a_{m-1}+d_{m-2}a_{m-2}+\ldots+d_0a_0$ (so that the associated string is 
$d_{m-1}d_{m-2}\ldots d_0$),
while, following the same scheme, the remaining integers (i.e. 
the integers $0\leq j<a_{m-1}$) have a representation with 
less than $m$ digits. For example: the representation of 
$a_{m-1}-1=d_{m-2}a_{m-2}+\ldots+d_0a_0$ has $m-1$ digits.
For our purpose, we require that all the 
representations of the considered integers $\ell\in\{0,1,2,\ldots,a_m-1\}$ have $m$ 
digits, so we pad the string on the left with 0's until we have $m$ digits: the 
representation of $a_{m-1}-1$ becomes 
$a_{m-1}-1=0a_{m-1}+d_{m-2}a_{m-2}+\ldots+d_0a_0$ (therefore, the associated string 
is $0d_{m-2}\ldots d_0$).

\medskip
With this little adjustment, we now define the following sets:

\vspace{.2cm}
$\mathscr L_0=\{\varepsilon\}$ ,

\vspace{.2cm}

$\mathscr L_m=\{d_{m-1}\ldots \ d_0\mid \mbox{the string}
\ d_{m-1}\ldots \ d_0 \ \mbox{is the representation of an integer}\\
\ell<a_m\ \mbox{in the numeration system } \{a_n\}_{n\geq 0}\ \}$.

\noindent
Finally, we denote by $\mathscr L$ the language obtained by taking the union of all the sets $\mathscr L_m$:

\vspace{.2cm}
$\mathscr L=\bigcup_{m\geq 0} \mathscr L_m$.

\bigskip
\noindent
We remark that each element of $\mathscr L_m$ has precisely $m$ 
digits, so that some string $d_{m-1}\ldots \ d_0$ can have a 
prefix consisting of consecutive zeros. 
Moreover, denoting by $|A|$ the cardinality of a set $A$,
it is $|\mathscr L_m|=a_m$. The $a_m$ elements are the representations of each $\ell \in\{0,1,\ldots, 
a_m-1\}$.  

\medskip
Referring to the sequence of Pell numbers $p_m=\{1,2,5,12,29,\ldots\}$ defined in Section 2, we have:
$$
\begin{array}{lll}
	
	\mathscr L_0&=&\{\varepsilon\} \\
	&&\\
	\mathscr L_1&=&\{0,1\} \\ 
	&&\\
	\mathscr L_2&=&\{00,01,10,11,20\}\\
	&&\\
	\mathscr L_3&=&\{000,001,010,011,020,100,101,
	110,111,120,200,201\}\\
	&&\\
	\mathscr L_4&=&\{0000,0001,0010,0011,0020,0100,0101,0110,0111,
	0120,0200,0201,1000,\\
	&&1001,1010,1011,1020,1100,1101,1110,1111,
	1120,1200,1201,2000,\\
	&&2001,2010,2011,2020\}\\
\end{array}
$$

The strings in $\mathscr L_2$ are, respectively, the representations of the integers $\ell\in \{0,1,2,3,4\}$. This corresponds to the case $m = 2$ where $a_m = 5$. Note that $\mathscr L_2$ contains exactly $a_2=5$ elements.

\bigskip
It is not difficult to realize that the alphabet of the language $\mathscr L$ strictly depends on the sequence $\{a_m\}_{m\geq 0}$. In general it is possible to set an upper bound for the digits $d_i$. From (\ref{minore}), we deduce $d_ia_i<a_{i+1}-\sum_{j=0}^{i-1} d_ja_j$, so that, since the numbers are all integers:

$$
d_ia_i\leq a_{i+1}-1-\sum_{j=0}^{i-1} d_ja_j\leq a_{i+1}-1\ ,
$$

\noindent leading to

\begin{equation}\label{bound}
	d_i\leq \left\lfloor \frac{a_{i+1}-1}{a_i} \right\rfloor\ .
\end{equation}

Therefore, the alphabet for $\mathscr L_m$ is given by
$\{0,1,\ldots ,s\} $ with
$$
s=\displaystyle \max_{i=0,1,\ldots, m-1}\left\{\left\lfloor \frac{a_{i+1}-1}{a_i} \right\rfloor\ \right\}\ ,
$$

\noindent
and, denoting by $\Sigma$ the alphabet for $\mathscr L$, we have $\Sigma=\{0,1,\ldots,t\}$ with

$$
t=\max_i\left\{\left\lfloor \frac{a_{i+1}-1}{a_i} \right\rfloor\ \right\}\ .
$$

\bigskip
\bigskip
In some recent papers (\cite{BBP1}, \cite{BBP2}) particular number sequences have been studied. More specifically, in \cite{BBP1} the recurrence

\begin{equation}
a_m=\left\{
\begin{array}{ll}
1& \mbox{if}\ m=0\\
k& \mbox{if}\ m=1\\
ka_{m-1}+ha_{m-2}& \mbox{if}\ m\geq 2\\
\end{array}\right.
\end{equation}

\noindent
with $k>h>0$ is considered. Here, the alphabet for the strings of $\mathscr L$ is $\{0,1,\ldots,k\}$ and it is possible to define a Gray code for them.

\medskip
In \cite{BBP2}, the following two-termed recurrence is analysed:

\begin{equation}\label{two_termed}
a_m=\left\{
\begin{array}{ll}
1& \mbox{if}\ m=0\\
k& \mbox{if}\ m=1\\
ka_{m-1}-ha_{m-2}& \mbox{if}\ m\geq 2\\
\end{array}\right.
\end{equation}

\noindent
with $k>h>0$. In this case, the alphabet is $\{0,1,\ldots,k-1\}$ and if $h$ is even it is possible to define a Gray code for the strings of $\mathscr L$. Moreover, recurrence (\ref{two_termed}) is equivalent to the following full history recurrence which involves only non-negative terms:

\begin{equation}\label{ful_history}
a_m=\left \{
\begin{array}{ll}
1 & \mbox{if}\ m=0\\
\\
k & \mbox{if}\ m=1\\
\\
(k-1)a_{m-1}+(k-h-1)a_{m-2}+\ldots&\vspace{.1cm}\\
+(k-h-1)a_1+(k-h-1)a_0+1& \mbox{if}\ m\geq 2\\
\end{array}\right.\ \ \ .
\end{equation}

\noindent
This is true also in the case $k=2$ and $h=0$ giving the one-termed recurrence

\begin{equation*}
a_m=\left\{
\begin{array}{ll}
1& \mbox{if}\ m=0\\
2a_{m-1}& \mbox{if}\ m\geq 1\\
\end{array}\right.\ \ \ ,
\end{equation*}

\noindent
equivalent to the full history recurrence
\begin{equation}\label{potenze}
a_m=\left \{
\begin{array}{ll}
1 & \mbox{if}\ m=0\\
a_{m-1}+a_{m-2}+\ldots+a_1+a_0+1& \mbox{if}\ m\geq 1\\
\end{array}\right.
\end{equation}

\noindent
(whose general term is $a_m=2^{m}$) , leading to the language 
$\mathscr L=\bigcup_{m\geq 0} \mathscr L_m$ where
$\mathscr L_m$ contains binary strings of length $m$. Denoting by $B_m$ the set of binary strings of length $m$ we have $\mathscr L_m\subseteq B_m$, but since $|\mathscr L_m|=|B_m|=2^m$, the two sets $\mathscr L_m$ and $B_m$ coincide for each $m\geq 0$.

\medskip
Let us introduce some notations (as in \cite{BBPSV}) useful throughout the rest of the paper.

\begin{itemize}
	\item If $\alpha$ is a symbol and $L$ is a list or a set of strings
	$L=(v_1,v_2,\ldots,v_s)$ or $L=\{v_1,v_2,\ldots,v_s\}$, then
	$\alpha\cdot L=(\alpha v_1,\alpha v_2\ldots,\alpha v_s)$ (or $\alpha\cdot L=\{\alpha v_1,\alpha v_2\ldots,\alpha v_s\}$) is the list or the set obtained by 
	left concatenating $\alpha$ to each string of $L$;
	\item if $i$ and $j$ are symbols, then $ij\cdot L$ is the list or the set obtained by left concatenating $i$ to each string of $j\cdot L$ (or equivalently $ij\cdot L = i \cdot (j \cdot L)$);
	\item if $L$ is a list or a set of strings, $\bar{L}$ is the list in the reverse order;
	\item if $L$ and $M$ are two lists, $L\circ M$ is their concatenation. For example, if $L=(v_1,v_2)$ and $M=(w_1,w_2)$, then $L\circ M=(v_1,v_2,w_1,w_2)$;
	\item if $L$ is a list, then $first(L)$ is the first element of $L$ and $last(L)$ is the last element of $L$.
\end{itemize}

It is known that the set $B_m$ can be defined by:

\begin{equation}\label{binario}
B_m=\left\{
\begin{array}{cl}
\{\varepsilon\},&\text{for } m=0\\
&\\
0\cdot B_{m-1}\cup
1\cdot B_{m-1}, & \text{for } m>0\\
\end{array}
\right.\ ,
\end{equation}

or equivalently by:

\begin{equation}\label{binario_full_history}
\footnotesize
B_m=\left\{
\begin{array}{cl}
\{\varepsilon\},&\text{for } m=0\\
&\\
0\cdot B_{m-1}\cup
10\cdot B_{m-2}\cup
\ldots
1^{m-2}0\cdot B_1\cup
1^{m-1}0\cdot B_0\cup
1^{m}, & \text{for } m>0\\
\end{array}\ \ ,
\right.
\end{equation}

where $B_i=\emptyset$ if $i<0$.

\medskip
A particular Gray code for the strings of $B_m$ is the \emph{local reflected lexicographical order} (\emph{lrl-order}, see \cite{C,V}), derived from the well known Binary Reflected Gray code \cite{G}. It is defined as:

\begin{equation}\label{BRGC}
\mathcal C_m=\left\{
\begin{array}{cl}
\varepsilon,&\text{for } m=0\\
&\\
0\cdot\bar{\mathcal C}_{m-1}\circ
1\cdot\mathcal C_{m-1}, & \text{for } m>0\\
\end{array}
\right.\ .
\end{equation}

\bigskip
\noindent
It is not difficult to see that $\mathcal C_m$ can be alternatively defined by:

\begin{equation}\label{BRGC_full_hystory}
\footnotesize
\mathcal C_m=\left\{
\begin{array}{cl}
\varepsilon,&\text{for } m=0\\
&\\
0\cdot\bar{\mathcal C}_{m-1}\circ
10\cdot\bar{\mathcal C}_{m-2}\circ
\ldots
1^{m-2}0\cdot\bar{\mathcal C}_1\circ
1^{m-1}0\cdot\bar{\mathcal C}_0\circ
1^{m}, & \text{for } m>0\\
\end{array}\ \ 
\right.
\end{equation}

where $\mathcal C_i=\emptyset$ if $i<0$.
\section{The case of $k$-genarilezed Fibonacci sequences}
We now investigate on the strings derived from the well known integer sequence of $k$-generalized Fibonacci numbers ($k\geq 2$). The sequence can be defined as follows:
\begin{equation}\label{kbonacci}
	f_{\ell}^{(k)}=\left \{
	\begin{array}{ll}
		2^{\ell} & \mbox{if} \ 0 \leq {\ell} \leq k-1\\
		\\
		f_{{\ell}-1}^{(k)}+f_{{\ell}-2}^{(k)}+\ldots+
		f_{{\ell}-k}^{(k)} & \mbox{if} \ {\ell} \geq k\ .
	\end{array}\right.
\end{equation}

\noindent
We observe that this recurrence can be seen as a particular case of (\ref{ful_history}) or (\ref{potenze}) (with suitable initial conditions) where the general term is given by the sum of the $k$ preceding terms rather than all the terms up to the first one.

Note that usually $k$-generalized Fibonacci sequences are defined with different initial conditions with respect to the ones we imposed. Here, in order to agree with the hypothesis of Theorem  \ref{teo}, we need that all the terms of the sequence are different, and this can be obtained by the above initial conditions.  

\medskip
We indicate the set of strings of length $m$ arising from the numeration system 
$\{f_{\ell}^{(k)}\}_{\ell\geq 0}$ by $\mathscr L_m^{(k)}$ and 
by $\mathscr L^{(k)}$ the language
$\mathscr L^{(k)}=\bigcup_{m\geq 0} \mathscr L_m^{(k)}$.

\bigskip
We have the following proposition:
\begin{proposition}\label{alphabet}
	The alphabet $\Sigma$ for the strings in $\mathscr L^{(k)}$ in the numeration system $\{f_n^{(k)}\}_{n\geq 0}$ is $\Sigma=\{0,1\}$.
	Moreover, the strings avoid the pattern $1^k$ (i.e. each string does not contain $k$ consecutive $1$'s).
\end{proposition}

\emph{Proof.}
From (\ref{bound}) we obtain
\begin{align*}
	d_i & \leq \left\lfloor \frac{f^{(k)}_{i+1}-1}{f^{(k)}_i}
	 \right\rfloor\ =
	 \left\lfloor\
	 \frac{f^{(k)}_{i}+f^{(k)}_{i-1}+\ldots +f^{(k)}_{i-k+1}-1}{f^{(k)}_i}
	 \right\rfloor\\
	 & = \left\lfloor 1 + \frac{f^{(k)}_{i-1}+\ldots + {f^{(k)}_{i-k+1}}-1}{f^{(k)}_i}
	 \right\rfloor\ .
\end{align*}

Since $\frac{f^{(k)}_{i-1}+\ldots +{f^{(k)}_{i-k+1}}-1}{f^{(k)}_i}<1$, we deduce that $\Sigma=\{0,1\}$.

\bigskip

Let $N$ be an integer and let $d_{m-1}d_{m-2}\ldots d_1d_0$ be its 
representation in the numeration system $\{f_n^{k}\}_{n\geq 0}$. We prove that it avoids the pattern $1^k$.
Suppose ad absurdum that $d_{m-1}d_{m-2}\ldots d_1d_0$ contains the 
pattern $1^k$ and let $j$ be the first index from the left such that 
$d_j=d_{j-1}=\ldots=d_{j-k+1}=1$ (clearly $k-1\leq j \leq m-1$
and, if $j<m-1$, it must be $d_{j+1}=0$).


From our assumption we have
\begin{align*}
N =
&\ d_{m-1}f^{(k)}_{m-1}+d_{m-2}f^{(k)}_{m-2}+\ldots
+d_{j+1}f^{(k)}_{j+1}\\
&+f^{(k)}_{j}+f^{(k)}_{j-1}+\ldots+f^{(k)}_{j-k+1}\\
&+d_{j-k}f^{(k)}_{j-k}+\ldots+d_1f^{(k)}_{1}+d_0f^{(k)}_{0}\ .\\
\end{align*}

Since $d_{j+1}=0$ and $f^{(k)}_{j}+f^{(k)}_{j-1}+\ldots+f^{(k)}_{j-k+1}=
f^{(k)}_{j+1}$, we obtain

\begin{align*}
&f^{(k)}_{j}+f^{(k)}_{j-1}+\ldots+f^{(k)}_{j-k+1}\\
&+d_{j-k}f^{(k)}_{j-k}+\ldots+d_1f^{(k)}_{1}+d_0f^{(k)}_{0}\geq f_{j+1}^{(k)}\\
\end{align*}

against Theorem \ref{teo} which assures 

\begin{align*}
&d_jf^{(k)}_{j}+d_{j-1}f^{(k)}_{j-1}+\ldots+d_{j-k+1}f^{(k)}_{j-k+1}\\
&+d_{j-k}f^{(k)}_{j-k}+\ldots+d_1f^{(k)}_{1}+d_0f^{(k)}_{0}< f_{j+1}^{(k)}\\
\end{align*}

Therefore the string $d_{m-1}d_{m-2}\ldots d_1d_0$ avoids the pattern $1^k$.

\cvd

\noindent
Denoting by $F_n^{(k)}$ the set of binary strings of length $n$ 
avoiding the pattern $1^k$, it is known \cite{K} that
$|F_n^{(k)}|=f_n^{(k)}$. From Proposition \ref{alphabet} we deduce 
that $\mathscr L_m^{(k)}\subseteq F_m^{(k)}$. Clearly, it is also
$|\mathscr L_n^{(k)}|=f_n^{(k)}$. Hence, the sets $\mathscr L_m^{(k)}$ and $F_m^{(k)}$ coincide. 

\bigskip

The sets $\mathscr L^{(k)}_m$ (and so $F_m^{(k)}$) can be defined recursively as follow:

\begin{equation}\label{linguaggio}
\mathscr L^{(k)}_m=\left\{
\begin{array}{cl}
B_m, &\text{for } m<k\\
&\\
0\cdot\mathscr L_{m-1}^{(k)}\cup
10\cdot\mathscr L_{m-2}^{(k)}\cup\ldots\cup
1^{k-1}0\cdot\mathscr L_{m-k}^{(k)}, & \text{for } m\geq k\\
\end{array}\ \ .
\right.
\end{equation}

The strings of $\mathscr L^{(k)}_m$ can be rearranged in a Gray Code
$\mathcal L^{(k)}_m$ with the Hamming distance equal to one (see \cite{V}):

\begin{equation}\label{vaino}
\mathcal L^{(k)}_m=\left\{
\begin{array}{cl}
\mathcal C_m,&\text{for } 0\leq m<k\\
&\\
0\cdot\bar{\mathcal L}_{m-1}^{(k)}\circ
10\cdot\bar{\mathcal L}_{m-2}^{(k)}\circ\ldots\circ
1^{k-1}0\cdot\bar{\mathcal L}_{m-k}^{(k)}, & \text{for } m\geq k\\
\end{array}\ \ .
\right.
\end{equation}

\section{Permutations}
Given a permutation $\pi\in S_m$, by inversion array of $\pi$ we mean the array $v(\pi)=v_1v_2\ldots v_{m-1}$ of dimension $m-1$ whose $i$-th entry counts the number of entries of $\pi$ at the right hand side of $\pi_i$ which are smaller than $\pi_i$. Formally, we have:
\begin{definition}\label{inversion_array}
If $\pi$ is a permutation of length $m$, the array 
$v(\pi)=v_1v_2\ldots v_{m-1}$ is the 
\emph{inversion array} of $\pi$, where
$$
v_i=\left|\{\pi_j|\pi_j<\pi_i, j>i\}\right|\ \mbox{for} \ i=1,2,\ldots, m-1\ .
$$
\end{definition}

\bigskip
We now associate the permutations in $S(321,312,23\ldots(k+1)1)$ (with $k\geq 2$) with the strings $\mathscr L^{(k)}$.
From Proposition \ref{alphabet} we know that a string $u\in \mathscr L^{(k)}$ avoids the pattern $1^k$.
We have the following proposition.

\begin{proposition}
Let $\pi$ be a permutation of length $m$ and let $v(\pi)$ the inversion array of $\pi$.
Then $\pi\in S_m(321,312,23\ldots(k+1)1)$ if and only if $v(\pi)\in \mathscr L^{(k)}_{m-1}$.
\end{proposition}

\emph{Proof.} We suppose $\pi\in S_m(321,312,23\ldots(k+1)1)$.
Since $\pi$ has length $m$, clearly $v(\pi)$ has length $m-1$, according to Definition \ref{inversion_array}. Since $\pi$ avoids $321$ and $312$, then for each entry $\pi_i$ there is at most only one entry $\pi_j$ with $j>i$ such that $\pi_i>\pi_j$, otherwise a pattern $321$ or $312$ would occur. Then, either $v_i=0$ or $v_i=1$ (therefore the alphabet of $v(\pi)$ is $\{0,1\}$). We now have to prove that $v(\pi)$ avoids $1^k$.

Let us suppose ad absurdum that such a pattern occurs and  $v_jv_{j+1}\ldots v_{j+k-1}$ be the leftmost occurrence of $1^k$ in $v(\pi)$ (it is $1\leq j\leq m-k$ and $v_{j-1}=0$ if $j\geq 1)$. 
Since $v_j=1$, there exists an index $r$ such that $\pi_{j+r}<\pi_j$. It must be $r\geq k$, otherwise, being $v_{j+r}=1$, it is $\pi_{j+r}>\pi_p$ for some $p>j+r$. But this is not possible since in this case a pattern $321$ would occur in $\pi_j\pi_{j+r}\pi_p$ (and the value of $v_j$ should be al least 2). Therefore, it is $r\geq k$.

Moreover, for each entry $\pi_{j+i}$ for $i=1,2,\ldots,k-1$ it is $\pi_{j+i}>\pi_{j+i-1}$ otherwise a pattern $321$ would occur in $\pi_{j+i-1}\pi_{j+i}\pi_{j+r}$ (and $v_{j+i-1}=2$ since $\pi_{j+i-1}>\pi_{j+i}$ and $\pi_{j+i-1}>\pi_{j+r}$). This implies that $\pi_j<\pi_{j+1}\ldots<\pi_{j+k-1}$. These entries of $\pi$ together with $\pi_{j+r}$ are an occurrence of $23\ldots (k+1)1$ against the hypothesis. Therefore $v(\pi)$ avoids $1^k$. 
 
\bigskip
Suppose now that $v(\pi)\in\mathscr L^{(k)}_{m-1}$. It means that either $v_i=0$ or $v_i=1$, then for any $\pi_i$ two entries $\pi_p$ and $\pi_q$ smaller than $\pi_i$, with $p,q>i$, do not exist. Therefore, the permutation $\pi$ avoid both  $321$ and $312$. We have to prove that $\pi$ avoids $23\ldots(k+1)1$, too.

Suppose ad absurdum that such a pattern occurs and let $\pi_{j_1}\pi_{j_2}\ldots\pi_{j_k}\pi_{j_{k+1}}$ be an occurrence of it. All the entries $\pi_p$ between $\pi_{j_i}$ and $\pi_{j_{i+1}}$ for $i=1,2,\ldots,k-1$ are such that $\pi_{j_i}<\pi_p<\pi_{j_{i+1}}$, otherwise either a pattern $312$ or $321$ occurs in $\pi_{j_{i}}\pi_p\pi_{j_{k+1}}$ if $\pi_p<\pi_{j_i}$, or a pattern $321$ occurs in $\pi_p\pi_{j_{i+1}}\pi_{j_{k+1}}$ if $\pi_p>\pi_{j_{k+1}}$.

This implies that there are at least $k$ 
consecutive entries 
$\pi_p\pi_{p+1}\ldots\pi_{p+k-1}$ of $\pi$ 
between $\pi_{j_i}$ and $\pi_{j_k}$ which are 
in increasing order 
($\pi_p<\pi_{p+1}<\ldots<\pi_{p+k-1}$) and for 
each of them it is $\pi_{p+r}>\pi_{j_{k+1}}$ 
for $r=0,1,\ldots, k-1$. Therefore in $v(\pi)$ 
we have $v_pv_{p+1}\ldots v_{p+k-1}=1^k$, 
against the hypothesis. Hence, $\pi$ avoids the pattern $23\ldots(k+1)1$, too. 

\cvd

Before giving the Gray code for $S_m(321,312,23\ldots(k+1)1)$, let us introduce some notations b.

Given a permutation $\pi\in S_n$, $\pi=\pi_1\pi_2\ldots \pi_n$, and a positive integer $p$,
we denote by $\pi\uparrow p=(\pi_1+p)(\pi_2+p)\ldots (\pi_n+p)$ the 
permutation of $[p+1,p+2,\ldots,p+n]$ obtained by $\pi$ by adding $p$ 
to each entry of $\pi$.

If $\rho\in S_p$, $\rho=\rho_1\rho_2\ldots\rho_p$ is a permutation of length $p$, then we denote by $\rho\cdot (\pi\uparrow p)=\rho_1\ldots\rho_p(\pi_1+p)\ldots(\pi_n+p)$ the permutation of
$[1,2,\ldots,n+p]$ obtained by concatenating $\rho$ with $\pi\uparrow p$.

If $\Pi=\{\pi^{(1)},\pi^{(2)},\ldots,\pi^{(\ell)}\}$ is a set of permutations, then the set of the  permutations $\pi^{(j)}\uparrow p$, for $j=1,2,\ldots ,\ell$ is denoted by $\Pi\uparrow p$.

\bigskip
We observe that the permutations of $S_m(321,312,23\ldots(k+1)1)$ can be recursively defined basing on Definition (\ref{linguaggio}) of $\mathscr L_{m-1}^{(k)}$. Indeed, if $v(\pi)\in \mathscr L_{m-1}^{(k)}$, then the inversion array $v(\pi')=0\cdot v(\pi)$ corresponds to the inversion array of the permutation $\pi'=1\cdot(\pi\uparrow 1)$ since a $0$ entry in an inversion array  corresponds to an entry of the permutation which is less than all the elements to its rigth. So, the first $0$ entry in $v(\pi')$ must correspond to the $1$ entry in $\pi'$.

More generally, by means of a similar argument, we can prove that if $v(\pi)\in \mathscr L_{m-j}^{(k)}$, then $v(\pi')=1^{j-1}0\cdot v(\pi)$ corresponds to the inversion array of the permutation $\pi'=23\ldots j1\cdot(\pi\uparrow j)$. Therefore, from (\ref{linguaggio}), by considering $B_m$ defined as in (\ref{binario_full_history}), we can define recursively the permutations of the set of $S_m(321,231,23\ldots(k+1)1)$ denoted, for short, by $S_m^{(k)}$. We give so the following definition:

\begin{equation}\label{permutazioni}
\footnotesize
S^{(k)}_m=\left\{
\begin{array}{cl}
\{\varepsilon\},&\text{for } m=0\\
&\\
1\cdot (S_{m-1}^{(k)}\uparrow 1)\cup
21\cdot (S_{m-2}^{(k)}\uparrow 2)\cup\ldots\cup
23\ldots k1\cdot (S_{m-k}^{(k)}\uparrow k), & \text{for } m\geq 1\\
\end{array}\ \ ,
\right.
\end{equation}
where $S_i^{(k)}=\emptyset$ if $i<0$.

\bigskip
We note that, given an inversion array $v(\pi)\in \mathscr L_{m-1}^{(k)}$ one can get the corresponding permutation $\pi\in S_m^{(k)}$ in a very simple way. Indeed, it is not difficult to prove that
denoting by $i_1,i_2,\ldots,i_r$ the indexes of the positions of the $0$'s in $v(\pi)\cdot 0$ (i.e. $v_{i_j}=0$ for $j=1,2,\ldots,r$), the entries of $\pi$ are:  

\begin{equation*}
\pi_i=\left\{
\begin{array}{cl}
1&\text{if } i=i_1\\
&\\
i_j+1 & \text{if } i=i_{j+1}\quad (j=1,2,\ldots,r-1)\\
&\\
i+1 & \text{if } i\neq i_j \quad (j=1,2,\ldots,r)\\
\end{array}\ \ .
\right.
\end{equation*}

\noindent
Moreover, if $\pi\in S_m^{(k)}$, the inversion array $v(\pi)=v_1v_2\ldots v_{m-1}$ is obtained as follows:
\begin{equation*}
v_i=\left\{
\begin{array}{cl}
0&\text{if } \pi_i\leq i\\
&\\
1 & \text{if } \pi_i>i\\
\end{array}\ \ \ .
\right.
\end{equation*}

\bigskip\noindent
Unlike in the case for strings, two permutations cannot differ
by a single position but, at least, by a transposition of two 
entries. One famous Gray code for $S_n$ (the set of the unrestricted 
permutations of length $n$) was given by Johnson \cite{J} and 
Trotter \cite{T}, where each permutation is obtained from the 
preceding one by a transposition of two consecutive entries.
 
As a final step of our paper we propose a Gray code $\mathcal S_m^{(k)}$ for the set  
$S_m^{(k)}$ which is induced by the 
Gray code $\mathcal L_m$ defined in (\ref{vaino}). We recursively define the list $\mathcal S_m^{(k)}$ as

\begin{equation}\label{permutazioni_Gray}
\footnotesize
\mathcal S^{(k)}_m=\left\{
\begin{array}{cl}
\varepsilon,&\text{for } m=0\\
&\\
1\cdot (\bar{\mathcal S}_{m-1}^{(k)}\uparrow 1)\circ
21\cdot (\bar{\mathcal S}_{m-2}^{(k)}\uparrow 2)\circ\ldots\circ
23\ldots k1\cdot (\bar{\mathcal S}_{m-k}^{(k)}\uparrow k), & \text{for } m\geq 1\\
\end{array}\ \ ,
\right.
\end{equation}
where $\mathcal S_i^{(k)}=\varepsilon$ if $i<0$.

\begin{proposition}
The list of permutation $\mathcal S_m^{(k)}$ is a Gray code where two consecutive permutations differ by a transposition of two consecutive entries. 
\end{proposition}

\emph{Proof.}
First of all we observe that, from (\ref{permutazioni_Gray}) for each $m\geq 2$, it is
\begin{equation}\label{qui}
\footnotesize
first(\mathcal S_m^{(k)})=
first(1\cdot (\bar{\mathcal S}_{m-1}^{(k)}\uparrow 1))=
1\cdot first(\bar{\mathcal S}_{m-1}^{(k)}\uparrow 1)=
1\cdot last(\mathcal S_{m-1}^{(k)}\uparrow 1)\ .
\end{equation}

We proceed by induction.
Easily, the lists $\mathcal S_0^{(k)}=\varepsilon$ and $\mathcal S_1^{(k)}=1$ and $\mathcal S_2^{(k)}=12,21$ are Gray codes where two consecutive permutations differ for a transposition of two consecutive entries (note that $\mathcal S_0^{(k)}$ and $\mathcal S_1^{(k)}$ are trivial cases).
Let us suppose that $\mathcal S_i^{(k)}$ is a Gray code
where two consecutive permutations differ by a transposition of two consecutive entries, for $i\leq m-1$.   

In order to prove that $\mathcal S_m^{(k)}$ is a Gray code 
we have to show only that the permutation
$last(23\ldots (j-1)1\cdot (\bar{\mathcal S}_{m-j+1}^{(k)}\uparrow (j-1)))$ 
and the permutation
$first(23\ldots j1\cdot (\bar{\mathcal S}_{m-j}^{(k)}\uparrow j))$
differ by a transposition of two consecutive entries, for $j=2,3,\ldots,k$.
Indeed, by the inductive hypothesis each 
$23\ldots j1\cdot (\bar{\mathcal S}_{m-j}^{(k)}\uparrow j)$ is a Gray code, for $j=1,2,\ldots,k$.

\medskip
We have:

\footnotesize
\begin{align*}
&last(23\ldots (j-1)1\cdot (\bar{\mathcal S}_{m-j+1}^{(k)}\uparrow (j-1))) =
23\ldots(j-1)1\cdot last((\bar{\mathcal S}_{m-j+1}^{(k)}\uparrow (j-1)))=\\
\\
&23\ldots(j-1)1\cdot first((\mathcal S_{m-j+1}^{(k)}\uparrow (j-1)))=
23\ldots(j-1)1\cdot (first(\mathcal S_{m-j+1}^{(k)})\uparrow (j-1))=\\
\\
&23\ldots(j-1)1\cdot ((1\cdot last(\mathcal S_{m-j}^{(k)}\uparrow 1))\uparrow (j-1)))=23\ldots (j-1)\bold {1j} \cdot last((\mathcal S_{m-j}^{(k)})\uparrow j)\ \ .
\\
\end{align*}

\normalsize\noindent
On the other hand:

\begin{align*}
&first(23\ldots (j-1)j1\cdot (\bar{\mathcal S}_{m-j}^{(k)}\uparrow j)) =
23\ldots (j-1)\bold{j1}\cdot last(\mathcal S_{m-j}^{(k)}\uparrow j)\\
\end{align*}

\noindent
so that the two permutations differ only for the transposition of the 
two consecutive entries $1j$ in the first one which become $j1$ in the second one.
\cvd

\section{Conclusion}

As already mentioned in Section \ref{intro}, this research was started in \cite{BR}, where interesting properties of the languages arising from particular linear recurrences were given.
Later (\cite{BBP1,BBP2}), the study has been enriched by the possibility of listing the languages in a Gray code order (at least in many cases), depending on the parameters of the linear recurrences. 
Here, we extend the results by starting from $k$-generalized Fibonacci recurrences and by involving pattern avoiding permutations, transferring on them the Gray code obtained for the language of strings by means of the inversion arrays.

A further improvement towards this direction could be devoted to  different particular recurrences, clearly leading to different languages. An interesting result would be the definition of a Gray code for them, as well as for pattern avoiding permutations possibly associated to.




\end{document}